\documentclass[draft]{amsart}
\usepackage[utf8]{inputenc}
\usepackage{amsthm}
\usepackage{amsmath}
\usepackage{amssymb}
\usepackage{graphicx}
\newtheorem{theorem}{Theorem}[section]
\newtheorem{lemma}[theorem]{Lemma}

\theoremstyle{definition}

\newtheorem{example}[theorem]{Example}

\newtheorem{remark}[theorem]{Remark}

\numberwithin{equation}{section}

\newtheorem*{thma}{Theorem A}

\begin{document}

\title[Quantization for self-similar measures with overlaps]{Asymptotic order of the quantization error for a class of self-similar measures with overlaps}
\author{Sanguo Zhu}

\address{School of Mathematics and Physics, Jiangsu University
of Technology\\ Changzhou 213001, China.}
\email{sgzhu@jsut.edu.cn}


\begin{abstract}
Let $\{f_i\}_{i=1}^N$ be a set of equi-contractive similitudes on $\mathbb{R}^1$ satisfying the finite-type condition. We study the asymptotic quantization error for the self-similar measures $\mu$ associated with $\{f_i\}_{i=1}^N$ and a positive probability vector. With a verifiable assumption, we prove that the upper and lower quantization coefficient for $\mu$ are both bounded away from zero and infinity. This can be regarded as an extension of Graf and Luschgy's result on self-similar measures with the open set condition. Our result is applicable to a significant class of self-similar measures with overlaps, including Erd\"{o}s measure, the $3$-fold convolution of the classical Cantor measure and the self-similar measures on some $\lambda$-Cantor sets.
\end{abstract}
\maketitle

\section{Introduction}
The quantization problem for a probability measure $\nu$ on $\mathbb{R}^q$ consists in the discrete approximation of $\nu$ by probability measures of finite support in $L_r$-metrics. This problem has a deep background in information theory and some engineering technology (cf. \cite{GN:98}). We refer to Graf and Luschgy \cite{GL:00} for rigorous mathematical foundations of quantization theory.

In the past years, the quantization problem has been extensively studied for fractal measures (cf. \cite{GL:01,GL:02,GL:04,LM:02,PK:01,PK:04,Zhu:18,Zhu:23}). With certain separation condition for the corresponding iterated function system (IFS), the asymptotics of the quantization error for self-similar measures have been well studied by Graf and Luschgy (cf. \cite{GL:01,GL:02,GL:04} ).
Up to now, very little is known about the asymptotics of the quantization error for self-similar measures with overlaps.

In this note, we study the quantization problem for the self-similar measures associated with a class of equi-contractive IFSs satisfying the finite type condition. Based on D.-J. Feng's work in \cite{Feng:03,Feng:05}, we determine the exact convergence order of the quantization error for a class of self-similar measures with overlapping structure.
\subsection{Asymptotics of the quantization errors}
Let $r\in (0,\infty)$ be given. For every $k\geq 1$, let $\mathcal{D}_k:=\{\alpha\subset\mathbb{R}^q:1\leq {\rm
card}(\alpha)\leq k\}$, where ${\rm card}(A)$ denotes the cardinality of a set $A$. Let $d$ denote the Euclidean metric on $\mathbb{R}^q$. The $k$th quantization
error for $\nu$ of order $r$ can be defined by
\begin{eqnarray}\label{quanerror}
e_{k,r}(\nu)=\inf\limits_{\alpha\in\mathcal{D}_k}\bigg(\int d(x,\alpha)^{r}d\nu(x)\bigg)^{1/r}.
\end{eqnarray}
By \cite{GL:00}, $e_{k,r}(\nu)$ agrees with the minimum error in the approximation of $\nu$ by probability measures supported on at most $k$ points in  $L_r$-metrics. Let $C_{k,r}(\nu)$ denote the set of all $\alpha\in \mathcal{D}_k$ such that the infimum in (\ref{quanerror}) is attained. Let $|x|$ denote the Euclidean norm of $x$. By \cite{GL:00}, $C_{k,r}(\nu)$ is non-empty whenever $\int |x|^rd\nu(x)<\infty$.

The asymptotics of the quantization error for $\nu$ can be characterized by the
upper and lower quantization coefficients of order $r$:
\[
\overline{Q}_r^s(\nu):=\limsup_{k\to\infty}k^{\frac{r}{s}}e_{k,r}^r(\nu),\;\underline{Q}_r^s(\nu):=\liminf_{k\to\infty}k^{\frac{r}{s}}e_{k,r}^r(\nu);\;s\in (0,\infty).
\]
The upper (lower) quantization dimension for $\nu$ of order $r$ is the critical point at which the upper (lower) quantization
coefficient jumps from zero to infinity:
\[
\overline{D}_r(\nu)=\limsup_{k\to\infty}\frac{\log k}{-\log e_{k,r}(\nu)};\;\;\underline{D}_r(\nu)=\liminf_{k\to\infty}\frac{\log k}{-\log e_{k,r}(\nu)}.
\]
When $\overline{D}_r(\nu)$ and $\underline{D}_r(\nu)$ agree, we say that the quantization dimension for $\nu$ of order $r$ exists and denote the common value by $D_r(\nu)$.

Compared with $\overline{D}_r(\nu)$ and $\underline{D}_r(\nu)$, people are more concerned about the upper and lower quantization coefficient, because they provide us with the exact order of the $n$th quantization error when they are both positive and finite.

Let $(f_i)_{i=1}^N$ be a family of similitudes on $\mathbb{R}^q$.
By \cite{Hut:81}, there exists a unique non-empty compact set $E$ satisfying
$E=\bigcup_{i=1}^Nf_i(E)$.
The set $E$ is called the self-similar set determined by $(f_i)_{i=1}^N$. Given a positive probability vector $(p_i)_{i=1}^N$, there exists a unique Borel probability measure $\mu$ satisfying $\mu=\sum_{i=1}^Np_i\mu\circ f_i^{-1}$. This measure is called the self-similar measure associated with $(f_i)_{i=1}^N$ and $(p_i)_{i=1}^N$. We say that $(f_i)_{i=1}^N$ satisfies the open set condition (OSC), if there exists some bounded non-empty open set $U$, such that $f_i(U),1\leq i\leq N$, are disjoint subsets of $U$.

Assuming the OSC, Graf and Luschgy established complete results for the asymptotics of the quantization error for self-similar measures (\cite{GL:01,GL:02}). The main difficulty, in the absence of the OSC, lies in the fact, that the hereditary law of the measures over cylinder sets can hardly be well tracked, due to the overlaps.

A recent breakthrough by Kesseb\"{o}hmer et al identified the upper quantization dimension of an arbitrary compactly supported probability measure with its R\'{e}nyi dimension at a critical point \cite{KNZ:22}. This work, along with Peres and Solomyak's results on the $L^q$-spectrum (cf. \cite{PS:00}), implies that, \emph{the quantization dimension for every self-similar measure on $\mathbb{R}^q$ exists}. Combining the results in \cite{KNZ:22} and those in \cite{Feng:05b,LN:98,NX:19}, one can obtain explicit formulas for the quantization dimension for a large class of self-similar measures with overlaps. However, the work in \cite{KNZ:22} does not provide us with exact convergence order for the quantization error. Therefore, we need to examine the finiteness and positivity of the quantization coefficient.

\subsection{Equi-contractive IFS and finite-type condition}

Let $0<\rho<1$ and $N\geq 2$. In the present paper, we consider the following IFS on $\mathbb{R}^1$:
\begin{equation}\label{equiifs}
f_i(x)=\rho x+b_j,\;\;0=b_1<b_2<\ldots<b_N=1-\rho.
\end{equation}
We call $(f_i)_{i=1}^N$ an equi-contractive IFS, since the contraction ratios are identical.
Let $|B|$ denote the diameter of a set $B$. We denote by $E$ the self-similar set associated with the IFS in (\ref{equiifs}), we clearly have $|E|=1$. Up to some suitable rescaling, the assumptions $b_1=0$ and $b_N=1-\rho$ can be removed (cf. Example \ref{3fold}).

Following D.-J. Feng \cite{Feng:03}, we say that $(f_i)_{i=1}^N$ satisfies the \emph{finite-type condition} (FTC), if there exists a finite set $\Gamma$, such that for $n\geq 1$, and every pair $\sigma,\omega\in\Omega_n$,
\[
{\rm either}\;\;\rho^{-n}|f_\sigma(0)-f_\omega(0)|\in\Gamma\;\;{\rm or}\;\;\rho^{-n}|f_\sigma(0)-f_\omega(0)|>1.
\]
One may see Ngai and Wang \cite{NW:01} for the FTC in a more general setting. For related work on the IFS satisfying the FTC, we refer to \cite{Feng:03,Feng:05b,Feng:05,Hare:16,Hare:18,NW:01,RW:97,Wu:22}.

In the study of the $L^q$-spectrum for self-similar measures, Feng \cite{Feng:03} proposed a method of partitioning the set $[0,1]$ into non-overlapping intervals and established characterizations for the hereditary law of the measure $\mu$ over such intervals. Feng's method and results will enable us to determine the asymptotics of the quantization error for a significant class of self-similar measures with overlaps.
\subsection{Statement of the main result}
  We write
\[
\mathcal{A}:=\{1,2,\ldots N\},\;\mathcal{A}_n:=\mathcal{A}^n,\;n\geq 1;\;\;\mathcal{A}^*:=\bigcup_{n\geq 1}\mathcal{A}_n.
\]
Let $\theta$ denote the empty word and $\mathcal{A}_0:=\{\theta\}$. We define $|\sigma|:=k$ for $\sigma\in\mathcal{A}_k$. For $n\geq h\geq 1$ and $\sigma=\sigma_1\ldots\sigma_n\in\mathcal{A}_n$, we write $\sigma|_h:=\sigma_1\ldots\sigma_h$. Define
\begin{eqnarray*}
f_\sigma:=\left\{\begin{array}{ll}
f_{\sigma_1}\circ f_{\sigma_2}\circ\cdots\circ f_{\sigma_n},&{\rm if}\;n\geq 1\\
{\rm id_{\mathbb{R}^1}}&{\rm if}\;n=0
\end{array}\right..
\end{eqnarray*}

We need the total self-similarity which is proposed by Broomhead, Montaldi and Sidorov (see \cite{BMS:04}): $E$ is totally self-similar if $f_\omega(E)=f_\omega([0,1])\cap E$ for every $\omega\in\mathcal{A}^*$. One may see \cite{BMS:04,DKY:19} for some interesting results and remarks on the total self-similarity.
Now we are able to state our main result.
\begin{theorem}\label{FTC}
Let $(f_i)_{i=1}^N$ be as defined in (\ref{equiifs}) satisfying the FTC. Let $E$ denote the self-similar set determined by $(f_i)_{i=1}^N$ and $\mu$ the self-similar measure associated with $(f_i)_{i=1}^N$ and a positive probability vector $(p_i)_{i=1}^N$. Assume that $E$ is totally self-similar. Then
for $s_r=D_r(\mu)$, we have
\begin{equation}\label{FTCeq}
0<\underline{Q}_r^{s_r}(\mu)\leq\overline{Q}_r^{s_r}(\mu)<\infty.
\end{equation}
\end{theorem}
The proof for Theorem \ref{FTC} relies on Feng's work \cite{Feng:03,Feng:05} and some results of Feng and Lau in \cite{Fenglau:02}. The measures as studied in \cite{RS:23} will be treated as a particular case of re-scaled $\lambda$-Cantor measures.
\section{Preliminaries}
In this section, we review some terminologies and known results of Feng, which we will work with in the remainder of the paper. We refer to \cite{Feng:03,Feng:05} for more details.

\subsection{Net intervals and characteristic vectors}
For every $n\geq 0$, we write
\[
P_n:=\{f_\sigma(0):\sigma\in\mathcal{A}_n\}\cup\{f_\sigma(1):\sigma\in\mathcal{A}_n\}\;\;{\rm and}\;\;t_n:={\rm card}(P_n).
 \]
 Let $(h_i)_{i=1}^{t_n}$, be the enumeration of the elements of $P_n$ in the increasing order. Define
\begin{equation}\label{netinterval}
\mathcal{F}_n:=\{[h_i,h_{i+1}]:(h_i,h_{i+1})\cap E\neq\emptyset, 1\leq i\leq t_n-1\}.
\end{equation}
The intervals in $\mathcal{F}_n$ are called \emph{net intervals} of order $n$.

For $\Delta=[0,1]$, let
$\ell_0(\Delta):=1, V_0(\Delta):=(0)$ and $r_0(\Delta):=1$.
 Now for $n\geq 1$ and $\Delta=[a,b]\in\mathcal{F}_n$, let $\ell_n(\Delta):=\rho^{-n}(b-a)$. Define
\[
\Upsilon_n(\Delta):=\{\rho^{-n}(a-f_\sigma(0)):\sigma\in\mathcal{A}_n, f_\sigma(E)\cap (a,b)\neq\emptyset\}.
\]

Let $(a_i)_{i=1}^k$, be the enumeration of $\Upsilon_n(\Delta)$, in the increasing order. Define
$v_n(\Delta):=k$ and $V_n(\Delta):=(a_i)_{i=1}^k$.
  Let $\hat{\Delta}\in\mathcal{F}_{n-1}$ such that $\Delta\subset\hat{\Delta}$. We denote by $(\Delta^i)_{i=1}^l$, the enumeration of all the sub-net-intervals in $\mathcal{F}_n$ of $\hat{\Delta}$ with $\ell_n(\Delta^i)=\ell_n(\Delta)$ and $V_n(\Delta^i)=V_n(\Delta)$, in the increasing order. Let $r_n(\Delta):=j$ for which $\Delta^j=\Delta$. The \emph{characteristic vector} for $\Delta$ is then defined by
$\mathcal{C}_n(\Delta):=(\ell_n(\Delta),V_n(\Delta),r_n(\Delta))$. By \cite[Lemma 2.2]{Feng:03}, the set $\Omega:=\{\mathcal{C}_n(\Delta):\Delta\in\mathcal{F}_n,n\geq 0\}$ is finite, whenever the FTC is fulfilled.
One may regard $\alpha=\mathcal{C}_n(\Delta)$ as the type for $\Delta\in\mathcal{F}_n$. We sometimes simply write $V(\alpha)$ for $V_n(\Delta)$ and write $v(\alpha)$ for $v_n(\Delta)$, because it depends on the type $\alpha$ rather than $\Delta$ itself. We define $\Omega^*:=\bigcup_{k\geq1}\Omega^k$.

\subsection{Admissible words}
For $n\geq 1$ and $\Delta\in\mathcal{F}_n$, there exists a unique finite sequence $(\Delta^{(i)})_{i=0}^n$ such that $\Delta^{(n)}=\Delta,\Delta^{(i)}\in\mathcal{F}_i$ and $\Delta^{(i)}\subset\Delta^{(i-1)}$ for every $1\leq i\leq n$. The sequence $(\mathcal{C}_i(\Delta^{(i)}))_{i=0}^n$ is called the \emph{symbolic expression} for $\Delta$.

Let $\alpha\in\Omega$, and $\Delta\in\mathcal{F}_n$ with $n\geq 0$ and $\mathcal{C}_n(\Delta)=\alpha$. Let $(\Delta^i)_{i=1}^k$ be the enumeration of sub-net-intervals of $\Delta$ of order $n+1$ in the increasing order. For $1\leq j\leq k$, we write $\mathcal{C}_{n+1}(\Delta^j)=\alpha_j$. Let
$\xi:\Omega\to\Omega^*: \xi(\alpha):=\alpha_1\ldots\alpha_k$. Define
\begin{eqnarray*}
A_{\alpha,\beta}:=\left\{\begin{array}{ll}
1&{\rm if}\;\beta=\alpha_i\;\;{\rm for\;some}\;i\\
0&{\rm otherwise}
\end{array}\right.,\;\;\beta\in\Omega.
\end{eqnarray*}
A word $\gamma_1\ldots\gamma_n\in\Omega^*$ is \emph{admissible} if $A_{\gamma_h,\gamma_{h+1}}=1$ for every $1\leq h\leq n-1$.
\subsection{Measures of net intervals}
Let $\hat{\Delta}=[c,d]\in\mathcal{F}_{n-1}$ and $\Delta=[a,b]\in\mathcal{F}_n$ with $\Delta\subset\hat{\Delta}$ and $\alpha=\mathcal{C}_{n-1}(\hat{\Delta}), \beta=\mathcal{C}_n(\Delta)$. We write
\[
v_{n-1}(\hat{\Delta})=k, v_n(\Delta)=l; V_{n-1}(\hat{\Delta})=(c_j)_{j=1}^k, V_n(\Delta)=(a_i)_{i=1}^l.
\]
For $1\leq j\leq k, 1\leq i\leq l$, let (see \cite[Lemma 3.2]{Feng:03})
\begin{eqnarray}\label{tji}
w_{j,i}:&=&\left\{\begin{array}{ll}
p_h&{\rm if}\;c-\rho^{n-1}c_j+\rho^{n-1}b_h=a-\rho^na_i\;{\rm for\;some}\;h\\
0&{\rm otherwise}
\end{array}\right..\\t_{j,i}:&=&w_{j,i}\mu([a_i,a_i+\ell_n(\Delta)])(\mu([c_j,c_j+\ell_{n-1}(\hat{\Delta})]))^{-1}.\nonumber
\end{eqnarray}
Define $T(\alpha,\beta):=(t_{j,i})_{k\times l}$. Let $\|\cdot\|_1$ denotes the $l_1$-norm of a vector. Let $\gamma_0\gamma_1\ldots\gamma_n$ be the symbolic expression for $\Delta\in\mathcal{F}_n$. By \cite[Theorem 3.3]{Feng:03},
\begin{equation}\label{heredi}
\mu(\Delta)=\|T(\gamma_0,\gamma_1)T(\gamma_1,\gamma_2)\cdots T(\gamma_{n-1},\gamma_n)\|_1.
\end{equation}

\subsection{Some of Feng's results}

 By \cite{Feng:05}, there exists exactly one \emph{essential class} $\hat{\Omega}\subset\Omega$, such that (i) for every $\alpha\in\hat{\Omega}$, we have $\beta\in\hat{\Omega}$ when $\alpha\beta$ is admissible; (ii) for every $\alpha,\beta\in\hat{\Omega}$, there exist $k\geq 0$ and $\gamma\in\hat{\Omega}^k$, such that $\alpha\gamma\beta$ is admissible.

Let $\hat{\Omega}=:\{\eta_1,\ldots,\eta_s\}$. We select an integer $n_0\geq 1$ and a net interval $I_0\in\mathcal{F}_{n_0}$ with $\mathcal{C}_{n_0}(I_0)=\eta_1$. Assume that $\gamma_0\ldots\gamma_{n_0-1}\eta_1$ is the symbolic expression for $I_0$. Write $\Theta_0:=\gamma_0\ldots\gamma_{n_0-1}$. Then for $\Delta\in\mathcal{F}_{n_0+k}$ with $\Delta\subset I_0$, its symbolic expression is of the following form: $\Theta_0\eta_1\eta_{i_1}\eta_{i_2}\ldots\eta_{i_k}$.

We identify $\hat{\Omega}$ with $\{1,2,\ldots,s\}$ and write $\sigma_1\ldots\sigma_n$ for $\eta_{\sigma_1}\ldots\eta_{\sigma_n}$. We define
\begin{eqnarray*}
&&\mathcal{B}_k:=\{\sigma\in\hat{\Omega}^k: \sigma_1=1, A_{\sigma_i,\sigma_{i+1}}=1\;{\rm for\;all}\;1\leq i\leq k-1\};\\
&&\mathcal{B}^*:=\bigcup\limits_{k\geq1}\mathcal{B}_k,\;\mathcal{B}_\infty:=\{\sigma\in\hat{\Omega}^\mathbb{N}: \sigma_1=1, A_{\sigma_i,\sigma_{i+1}}=1\;{\rm for\;all}\;i\geq 1\}.
\end{eqnarray*}
For every $\sigma\in\mathcal{B}^*$, let $\Delta_\sigma$ denote the net interval with symbolic expression $\Theta_0\sigma$.

For $k\geq 1$, let $\textbf{e}$ denote $k$-dimensional column vector with all entries equal to $1$ and $\textbf{e}^T$ its transpose. For a $k\times k$ matrix $B$, let $\|B\|:=\textbf{e}^TB\textbf{e}$, For $\sigma\in\hat{\Omega}^n$, we write $B_\sigma$ for the product $B_{\sigma_1}\cdot B_{\sigma_2}\ldots \cdot B_{\sigma_n}$ of $k\times k$ matrices $B_{\sigma_i}, 1\leq i\leq n$.

\begin{thma}(See \cite[Proposition 5.1]{Feng:05}.)
Let $T:=\sum_{i=1}^sv(\eta_i)$. There exist non-negative $T\times T$ matrices $M_i,1\leq i\leq s$, such that
(1) for every $\sigma\in\hat{\Omega}^n$, $M_{\sigma}\neq 0$ if and only if $\sigma$ is admissible;
(2) $(M_i)_{i=1}^s$ is irreducible in the sense that there exists a positive integer $p$ such that $\sum_{k=1}^p(\sum_{i=1}^sM_i)^k>0$;
(3) there exist constants $C_1,C_2>0$, such that for every $\sigma\in\mathcal{B}_n$, we have
\begin{eqnarray*}\label{feng}
({\rm f1})\; C_1\|M_\sigma\|\leq \mu(\Delta_\sigma)\leq C_1^{-1}\|M_\sigma\|;\;({\rm f2})\;C_2\rho^n\leq|\Delta_\sigma|\leq C_2^{-1}\rho^n.
\end{eqnarray*}
\end{thma}
\subsection{Some remarks}
For $\sigma\in\mathcal{B}^*$, let $|\sigma|, \sigma|_h$ be defined in the same way as for the words in $\mathcal{A}^*$. For $\sigma,\omega\in\mathcal{B}^*$ with $|\sigma|\leq|\omega|$ and $\sigma=\omega|_{|\sigma|}$, we write $\sigma\prec\omega$. Define
\begin{eqnarray*}
&&\sigma^\flat:=\left\{\begin{array}{ll}
\theta,&{\rm if}\;|\sigma|=1\\
\sigma|_{|\sigma|-1}&{\rm if}\;|\sigma|>1
\end{array}\right.;\;\mathcal{E}_r(\sigma):=\mu(\Delta_\sigma)|\Delta_\sigma|^r,\;\sigma\in\mathcal{B}^*.
\end{eqnarray*}
 Let $A^\circ$ denote the interior of a set $A\subset\mathbb{R}^1$ and $A^c$ its complement.
\begin{remark}\label{descendant}
We have $\max_{\alpha\in\hat{\Omega}}|\xi(\alpha)|\geq 2$. In fact, by (\ref{netinterval}), we have, $I_0^\circ\cap E\neq\emptyset$. Thus, there exists an $\omega_0\in\mathcal{A}^*$ such that $f_{\omega_0}(E)\subset I_0^\circ$, implying ${\rm card}(I_0\cap E)=\infty$. Assume that $\max_{\alpha\in\hat{\Omega}}|\xi(\alpha)|=1$. By (f2), for $\sigma\in\mathcal{B}_\infty$, we have $|\Delta_{\sigma|_n}|\to 0$ as $n\to\infty$. Thus, the set $E\cap I_0$ would be a finite set, a contradiction. This can also be easily seen by considering different cases of the endpoints of $I_0$.
\end{remark}

\begin{remark}\label{z1}
(i) Assume that $T(j,k)\textbf{e}>0$ for every $jk\in\mathcal{B}_2$. Let $R_{j,k}^{(i)}$ denote the $i$th row of $T(j,k)$. Let
$C_3:=\min\{\|R_{j,k}^{(i)}\|_1: 1\leq i\leq v(j),jk\in\mathcal{B}_2\}$.
As an easy consequence of (\ref{heredi}), for every $\sigma\in\mathcal{B}^*$ with $|\sigma|\geq 2$, we have
$\mu(\Delta_\sigma)\geq C_3\mu(\Delta_{\sigma^\flat})$. This, \cite[Proposition 2.2]{Fenglau:09} and Remark \ref{descendant} further imply that $C_3\leq 2^{-1}$.
(ii) It was observed in \cite[p. 346]{Hare:16} that, when $E=[0,1]$, the assumption in (i) is fulfilled.
\end{remark}
\section{Proof of Theorem \ref{FTC}}

Let $(f_i)_{i=1}^N$, $E$ and $\mu$ be the same as in Theorem \ref{FTC}. Let $I_0$ be as selected in Section 2. We define $\mu_0:=\mu(\cdot|I_0)$ as the conditional measure of $\mu$ on $I_0$.
We will establish some estimates for the quantization error for $\mu_0$. By applying some auxiliary measures from \cite{Fenglau:02}, we will first prove (\ref{FTCeq}) for $\mu_0$, and then transfer this result to $\mu$ by applying the self-similarity of $\mu$.
Our first lemma shows that, when $E$ is totally self-similar, the assumption in Remark \ref{z1} (i) is fulfilled.
\begin{lemma}\label{g9}
Assume that $E$ is totally self-similar. Then for every pair $\alpha,\beta\in\Omega$ with $A_{\alpha,\beta}=1$, we have $T(\alpha,\beta)\textbf{e}>0$.
\end{lemma}
\begin{proof}
Let $\alpha,\beta\in\Omega$ with $A_{\alpha,\beta}=1$. We pick net intervals $\Delta=[a,b]\in\mathcal{F}_n$ and $\hat{\Delta}=[c,d]\in\mathcal{F}_{n-1}$ such that
$\Delta\subset\hat{\Delta}, \mathcal{C}_{n-1}(\hat{\Delta})=\alpha$ and $\mathcal{C}_n(\Delta)=\beta$.
We write $v(\alpha)=:p,\;v(\beta)=:l$, and $V(\alpha)=(c_j)_{j=1}^p,V(\beta)=(a_i)_{i=1}^l$.
For every $1\leq j\leq p$, there exists some $\hat{\sigma}\in\mathcal{A}_{n-1}$ such that
$c_j=\rho^{-(n-1)}(c-f_{\hat{\sigma}}(0))$ and $f_{\hat{\sigma}}(E)\cap (c,d)\neq\emptyset$.
By the definition of net intervals, we have,
$(a,b)\cap E\neq\emptyset$ and $[c,d]\subset f_{\hat{\sigma}}([0,1])$.
Using the total self-similarity of $E$, we deduce
\[
(a,b)\cap E\subset f_{\hat{\sigma}}([0,1])\cap E=f_{\hat{\sigma}}(E)=\bigcup_{k=1}^Nf_{\hat{\sigma}\ast k}(E).
\]
Thus, there exists some $1\leq k\leq N$, such that $(a,b)\cap f_{\hat{\sigma}\ast k}(E)\neq\emptyset$. It follows that for some $1\leq i\leq l$, we have
$f_{\hat{\sigma}\ast k}(0)=a-\rho^n a_i$. Hence,
$a-\rho^n a_i=c-\rho^{n-1}c_j+\rho^{n-1}b_k$.
By (\ref{tji}), we see that $w_{j,i}=p_k>0$. Hence, $T(\alpha,\beta)\textbf{e}>0$.
\end{proof}

In the following, we always assume that the assumption in Remark \ref{z1} (i) holds. For every $x\in\mathbb{R}$, let $[x]$ denote the largest integer not exceeding $x$. To obtain some estimates for the quantization error for $\mu_0$, we need the following lemma.
\begin{lemma}\label{z2}
Let $L\in\mathbb{N}$ and $B\subset\mathbb{R}^1$ with ${\rm card}(B)=L$. There exists a positive number $Z_{L,r}$, which is independent of $\sigma\in\mathcal{B}^*$, such that
\[
\int_{\Delta_\sigma}d(x,B)^rd\mu(x)\geq Z_{L,r} \mathcal{E}_r(\sigma).
\]
\end{lemma}
\begin{proof}
Let $k\geq 0,h\geq 1$, and $\sigma\in\mathcal{B}_{k+1}$ be given. We define
\begin{eqnarray*}
\Gamma(\sigma,h):=\{\omega\in\mathcal{B}^*:\sigma\prec\omega, |\omega|=|\sigma|+h\}.\label{g1}
\end{eqnarray*}
Since $\hat{\Omega}$ is an essential class, the  matrix $A=(A_{i,j})_{i,j=1}^s$ is irreducible. Using this and Remark \ref{descendant}, we deduce that, there exists some positive integer $H\leq s$, such that ${\rm card}(\Gamma(\sigma,H))\geq 2$. Note that every net interval has at least one sub-net-interval of the next order. Inductively, ${\rm card}(\Gamma(\sigma,h))\geq 2^{[h/s]}$ for every $h\geq s$. Let $k_L$ denote the smallest integer such that $2^{[k_L/s]}>3L+1$.
Because net intervals of the same order are pairwise non-overlapping, for every $b\in B$, we have
\[
{\rm card}(\{\tau\in\Gamma(\sigma,k_L): d(b,\Delta_\tau)\leq C_2\rho^{n_0+k+k_L}\})\leq 3.
\]
Using this and  Theorem A (f2), we may select some $\tau\in\Gamma(\sigma,k_L)$, such that
\[
d(\Delta_\tau,B)\geq C_2\rho^{-(n_0+k+k_L)}.
\]
This, together with Remark \ref{z1} and (f2), yields that
\begin{eqnarray*}
\int_{\Delta_\sigma}d(x,B)^rd\mu(x)\geq \int_{\Delta_\tau}d(x,B)^rd\mu(x)\geq C_3^{k_L}C_2^{2r}\rho^{k_Lr}\mathcal{E}_r(\sigma).
\end{eqnarray*}
Thus, the lemma is fulfilled with $Z_{L,r}:=C_3^{k_L}C_2^{2r}\rho^{k_Lr}$.
\end{proof}

Let $\eta_r:=C_3C_2^{2r}\rho^r$. For every $k\geq 1$, we define
\begin{eqnarray}\label{lambdakr}
\Lambda_{k,r}:=\{\sigma\in\mathcal{B}^*: \mathcal{E}_r({\sigma^\flat})\geq \eta_r^k\mu(I_0)I_0^r>\mathcal{E}_r(\sigma)\}.
\end{eqnarray}
Let $\phi_{k,r}:={\rm card}(\Lambda_{k,r})$.  By using Lemma \ref{z2} and \cite[Lemma 3]{KZ:16}, we are able to establish some estimates for the quantization error for $\mu_0$ of order $r$.
\begin{lemma}\label{ss1}
There exist constants $C_{4,r},C_5>0$, such that
\[
C_{4,r}\sum_{\sigma\in\Lambda_{k,r}}\mathcal{E}_r(\sigma)\leq e_{\phi_{k,r},r}^r(\mu_0)\leq C_5\sum_{\sigma\in\Lambda_{k,r}}\mathcal{E}_r(\sigma).
\]
\end{lemma}
\begin{proof}
For every $\sigma\in\Lambda_{k,r}$, let $a_\sigma$ be an arbitrary point in $\Delta_\sigma$. We have
\begin{eqnarray*}\label{g3}
e_{\phi_{k,r},r}^r(\mu_0)\leq\sum_{\sigma\in\Lambda_{k,r}}\int_{\Delta_\sigma} d(x,a_\sigma)^rd\mu_0(x)\leq \mu(I_0)^{-1}\sum_{\sigma\in\Lambda_{k,r}}\mathcal{E}_r(\sigma).
\end{eqnarray*}
As in Lemma \ref{z2}, for each $\sigma\in\Lambda_{k,r}$, we may choose a $\tau_\sigma\in\Gamma(\sigma,k_1)$ such that
\begin{eqnarray*}\label{s4}
d(\Delta_{\sigma\ast\tau_\sigma},\Delta_\sigma^c)\geq C_2|\Delta_{\sigma\ast\tau_\sigma}|;\;\mu(\Delta_{\sigma\ast\tau_\sigma})\geq C_3^{k_1}\mu(\Delta_\sigma)).
\end{eqnarray*}
Using this and (\ref{lambdakr}), for every pair $\sigma,\omega\in\Lambda_{k,r}$ of distinct words, we deduce
\begin{eqnarray}\label{s5}
&&\mathcal{E}_r(\sigma\ast\tau_\sigma)\geq C_3^{k_1}C_2^{2r}\rho^{k_1r}\eta_r\mathcal{E}_r(\omega\ast\tau_\omega);\\
&&d(\Delta_{\sigma\ast\tau_\sigma},\Delta_{\omega\ast\tau_\omega})\geq C_2\max(|\Delta_{\sigma\ast\tau_\sigma}|,|\Delta_{\omega\ast\tau_\omega}|).\label{s6}
\end{eqnarray}
Let $B_{k,r}:=\bigcup_{\sigma\in\Lambda_{k,r}}\Delta_{\sigma\ast\tau_\sigma}$. By Remark \ref{z1}, one gets
$\mu(B_{k,r})\geq C_3^{k_1}\mu(I_0)$. We define $\mu_{k,r}:=\mu(\cdot|B_{k,r})$.
By (\ref{g1})-(\ref{s6}), Lemma \ref{z2} and \cite[Lemma 3]{KZ:16},  there exists a constant $D=D(r)>0$, which is independent of $k$, such that
\begin{eqnarray}\label{g2}
e_{\phi_{k,r},r}^r(\mu_{k,r})\geq D\sum_{\sigma\in\Lambda_{k,r}}\mu_{k,r}(\Delta_{\sigma\ast\tau_\sigma})|\Delta_{\sigma\ast\tau_\sigma}|^r.
\end{eqnarray}
Let $\beta\in C_{\phi_{k,r},r}(\mu_0)$. Using (\ref{g2}), Theorem A (3) and Remark \ref{z1}, we deduce
\begin{eqnarray*}\label{t1}
e_{\phi_{k,r},r}^r(\mu_0)&\geq&\int_{B_{k,r}} d(x,\beta)^rd\mu_0(x)\\&=&\mu(I_0)^{-1}\mu(B_{k,r})\int_{B_{k,r}} d(x,\beta)^rd\mu_{k,r}(x)\\&\geq&  \mu(I_0)^{-1}\mu(B_{k,r})e_{\phi_{k,r},r}^r(\mu_{k,r})\\&\geq& DC_3^{k_1}C_2^{2r}\rho^{k_1r}\sum_{\sigma\in\Lambda_{k,r}}\mathcal{E}_r(\sigma).
\end{eqnarray*}
The lemma follows by defining $C_{4,r}:=DC_3^{k_1}C_2^{2r}\rho^{k_1r}$ and $C_5:=\mu(I_0)^{-1}$.
\end{proof}

Let $M_i,1\leq i\leq s$, be the same as in Theorem A. We write
$\widetilde{M}_{i,r}:=\rho^r M_i$.
Since the matrix norm is sub-multiplicative, in view of Theorem A (1), we define
\[
\Phi_r(t)=\lim_{n\to\infty}\frac{1}{n}\log\sum_{\sigma\in\hat{\Omega}^n} \|\widetilde{M}_{\sigma,r}\|^t,\;t>0.
\]
The function $\Phi_r$ corresponds to the pressure function $P$ as defined in \cite{Fenglau:02}.

For two variables $X,Y$ taking values in $(0,\infty)$, we write $X\asymp Y$ if there exists some constant $C>0$ such that $CY\leq X\leq C^{-1}Y$. The following facts are implied in the proof of \cite[Proposition 5.7]{Feng:03}:
\begin{eqnarray}\label{f1}
&&\sum_{\sigma\in\hat{\Omega}^n} \|\widetilde{M}_{\sigma,r}\|^t\asymp\sum_{\sigma\in\mathcal{B}_n} \|\widetilde{M}_{\sigma,r}\|^t,\label{f1}\\&&
\Phi_r(t)=\lim_{n\to\infty}\frac{1}{n}\log\sum_{\sigma\in\mathcal{B}_n} \|\widetilde{M}_{\sigma,r}\|^t,\;t>0.\label{f2}
\end{eqnarray}
 Next, using (\ref{f2}), we show that the function $\Phi_r$ has a unique zero (in $(0,1)$).
\begin{lemma}\label{z3}
There exists a unique $\xi_r\in(0,1)$ such that $\Phi(\xi_r)=0$. As a consequence, there exists a unique $s_r>0$ such that $\Phi_r(\frac{s_r}{s_r+r})=0$.
\end{lemma}
\begin{proof}
By Theorem A (f1), we have, $\sum_{\sigma\in\mathcal{B}_n} \|M_\sigma\|\asymp\mu(I_0)$. This implies that
$\Phi_r(1)\leq r\log\rho<0$. Since $N\geq 2$, we have
\[
\dim_B (I_0^\circ\cap E)=\dim_B E\geq \min\big(1,-\frac{\log2}{\log\rho}\big)=:2d_0>0.
\]
 Note that $I_0^\circ\cap E\subset\bigcup_{\sigma\in\mathcal{B}_n}\Delta_\sigma$ for $n\geq 1$.
Thus, using (f2), one can easily see that ${\rm card}(\mathcal{B}_n)>\rho^{-nd_0}$ for every large $n$. It follows that $\Phi_r(0)>0$.
Now let $\epsilon\in (0,1)$. By (f1) and Remark \ref{z1}, we deduce (cf. \cite[Lemma 5.2]{Fal:97})
\begin{eqnarray*}
\epsilon(\log C_3+r\log\rho)\leq\Phi_r(t+\epsilon)-\Phi_r(t)
 \leq r\epsilon\log\rho.
\end{eqnarray*}
Therefore, $\Phi_r(s)$ is strictly decreasing and continuous. The lemma follows.
\end{proof}
\begin{remark}
 Let $\tau(q)$ be as given in \cite[Proposition 5.7]{Feng:03}. The number $s_r$ agrees with $\frac{rq_r}{1-q_r}$, where $q_r$ is the unique number satisfying $-\tau(q_r)=rq_r$. By \cite[Theorem 1.11]{KNZ:22}, we have $D_r(\mu)=s_r$. This is independently implied by (\ref{FTCeq}).
\end{remark}
Inspired by \cite[Lemma 5.3]{Feng:03}, we are able to establish a relationship between the upper and lower quantization coefficients for $\mu$ and those for $\mu_0$. That is,

\begin{lemma}\label{z4}
There exists some $C_{6,r},C_7>0$, such that for all $t>0$, we have
\begin{eqnarray*}
&&C_7\overline{Q}_r^t(\mu_0)\leq \overline{Q}_r^t(\mu)\leq C_{6,r}\overline{Q}_r^t(\mu_0);\\
&&C_7\underline{Q}_r^t(\mu_0)\leq \underline{Q}_r^t(\mu)\leq C_{6,r}\underline{Q}_r^t(\mu_0).
\end{eqnarray*}
\end{lemma}
\begin{proof}
Let $\omega_0\in\mathcal{A}^*$, be the same as in Remark \ref{descendant}. Let $\alpha\in C_{n,r}(\mu)$. We have
\begin{eqnarray*}
e^r_{n,r}(\mu)\geq \int_{I_0} d(x,\alpha)^rd\mu(x)\geq\mu(I_0)e^r_{n,r}(\mu_0)\geq p_{\omega_0}e^r_{n,r}(\mu_0).
\end{eqnarray*}
Let $n\geq1$ and $B_n\in C_{n,r}(\mu_0)$. We define $\gamma_n:=f_{\omega_0}^{-1}(B_n)$. Then we have
\begin{eqnarray*}
e_{n,r}^r(\mu_0)&=&\int_{I_0} d(x,B_n)^rd\mu_0(x)\\&\geq&
\mu(I_0)^{-1}\sum_{\tau\in\mathcal{A}_{|\omega_0|}}p_\tau\int_{f_{\omega_0}(E)} d(x,B_n)^rd\mu\circ f_{\tau}^{-1}(x)
\\&\geq&\mu(I_0)^{-1}p_{\omega_0}\int_{f_{\omega_0}(E)} d(x,B_n)^rd\mu\circ f_{\omega_0}^{-1}(x)
\\&=&\mu(I_0)^{-1}p_{\omega_0}\rho^{|\omega_0|r}\int d(x,\gamma_n)^rd\mu(x)
\\&\geq&p_{\omega_0}\rho^{|\omega_0|r}e_{n,r}^r(\mu).
\end{eqnarray*}
It is sufficient to define $C_{6,r}:=p_{\omega_0}^{-1}\rho^{-|\omega_0|r}$ and $C_7:=p_{\omega_0}$.
\end{proof}

\emph{Proof of theorem \ref{FTC}}
Let $s_r$ be as defined in Lemma \ref{z3}. Let $(M_i)_{i=1}^s$ be the matrices in Theorem A. Since $(M_i)_{i=1}^s$ is irreducible, so is $(\widetilde{M}_{i,r})_{i=1}^s$. For every $n\geq 1$ and $\sigma\in\hat{\Omega}^n$, we define
$[\sigma]:=\{\tau\in\hat{\Omega}^{\mathbb{N}}:\tau|_n=\sigma\}$.
In view of Theorem A (1), we apply \cite[Theorem 3.2]{Fenglau:02} and deduce that, there exists a Borel probability measure $W$ on $\hat{\Omega}^{\mathbb{N}}$ such that, for every $n\geq1$, we have
\begin{equation}\label{t3}
W([\sigma])\asymp\|\widetilde{M}_{\sigma,r}\|^{\frac{s_r}{s_r+r}},\;\sigma\in\hat{\Omega}^n.
\end{equation}
Hence, $W([\sigma])\asymp\mathcal{E}_r(\sigma)^\frac{s_r}{s_r+r}$ for $\sigma\in\mathcal{B}^*$. Using this, (\ref{f1}) and (\ref{t3}), we deduce
\begin{eqnarray}
\sum_{\sigma\in\Lambda_{k,r}}(\mathcal{E}_r(\sigma))^{\frac{s_r}{s_r+r}}&\asymp&\sum_{\sigma\in\Lambda_{k,r}}W([\sigma])=
\sum_{\sigma\in\mathcal{B}_k}W([\sigma]) \asymp \sum_{\sigma\in\mathcal{B}_k}
\|\widetilde{M}_{\sigma,r}\|^{\frac{s_r}{s_r+r}}\nonumber\\&\asymp&\sum_{\sigma\in\hat{\Omega}_k}
\|\widetilde{M}_{\sigma,r}\|^{\frac{s_r}{s_r+r}}\asymp\sum_{\sigma\in\hat{\Omega}_k}W([\sigma])=1.\label{t2}
\end{eqnarray}
On the other hand, using Remark \ref{descendant} and along the line in \cite[Lemma 2.4]{Zhu:13}, one can check that $\phi_{k,r}\asymp\phi_{k+1,r}$.
Thus, by Lemma \ref{ss1}, (\ref{t2}) and \cite[Lemma 3.4]{Zhu:23}, we obtain that $0<\underline{Q}_r^{s_r}(\mu_0)\leq\overline{Q}_r^{s_r}(\mu_0)<\infty$.
Combining this and Lemma \ref{z4}, we obtain Theorem \ref{FTC}.

\section{Examples}
\begin{example}
Let $\rho=\frac{1}{2}(\sqrt{5}-1)$. Erd\"{o}s measure is the distribution measure of the random variable $(1-\rho)\sum_{n=0}^\infty\rho^n X_n$, where $X_n,n\geq 0$, are i.i.d random variables taking values $0$ and $1$ with probability $\frac{1}{2}$. This measure is exactly the self-similar measure associated with $(\frac{1}{2},\frac{1}{2})$ and the IFS:
$f_1(x)=\rho x, f_2(x)=\rho x+1-\rho$ (cf. \cite{LN:98}). Since $E=[0,1]$ is totally self-similar, (\ref{FTCeq}) holds by Theorem \ref{FTC}.
\end{example}

\begin{example}\label{3fold}
The Cantor measure $\zeta$ is the self-similar measure associated with the probability vector $(\frac{1}{2},\frac{1}{2})$ and the following IFS:
$f_1(x)=\frac{1}{3}x,\;f_2(x)=\frac{1}{3}x+\frac{2}{3}$.
As is noted in \cite{FLN:02}, the $3$-fold convolution $\mu=\zeta\ast\zeta\ast\zeta$ agrees with the self-similar measure associated with $\textbf{P}=(\frac{1}{8},\frac{3}{8},\frac{3}{8},\frac{1}{8})$ and the IFS:
$g_i(x)=\frac{x}{3}+\frac{2}{3}i,\;i=0,\;1,\;2,\;3$.
We have, (\ref{FTCeq}) holds. To see this, we define
\[
\varphi(x):=\frac{1}{3}x;\;\;h_i(x)=\frac{x}{3}+\frac{2}{9}i,\;\;i=0,\;1,\;2,\;3;\;x\in\mathbb{R}.
\]
By \cite[Theorem 2.9]{NW:01}, $(h_i)_{i=0}^3$ satisfies the FTC.
Let $\nu$ denote the self-similar measure associated with $(h_i)_{i=0}^3$ and $\textbf{P}$.
Note that $\varphi\circ g_i(x)=h_i\circ\varphi(x)$ for $i=0,1,2,3$. By induction, we obtain
$\varphi\circ g_\sigma(x)=h_\sigma\circ\varphi(x)$ for every $\sigma\in\mathcal{A}^*$, where $\mathcal{A}:=\{0,1,2,3\}$.
Using this, one can check that
 $\mu=\nu\circ\varphi$. Note that ${\rm supp}(\nu)=E=[0,1]$ and $\varphi$ is a similitude. By Theorem \ref{FTC}, (\ref{FTCeq}) holds.
\end{example}

\begin{example}\label{lambdacantor}
Let $\lambda\in (0,1)$. We consider the following IFS:
\begin{equation}\label{lambdacantoreq}
f_1(x)=\frac{1}{3}x,\;\;f_2(x)=\frac{1}{3}x+\frac{\lambda}{3},\;\;f_3(x)=\frac{1}{3}x+\frac{2}{3}.
\end{equation}
The self-similar set $E_\lambda$ associated with the above IFS is called a $\lambda$-Cantor set (cf. \cite{RW:97}). Let $\mu$ denote the self-similar measure associated with $(f_i)_{i=1}^3$ and a positive probability vector $(p_i)_{i=1}^3$. By \cite[Theorem 1]{DKY:19}, $E_\lambda$ is totally self-similar if and only if $\lambda=1-3^{-m}$ for some $m\in\mathbb{N}$. For such a $\lambda$, the IFS in (\ref{lambdacantoreq}) clearly satisfies the FTC by \cite[Theorem 2.9]{NW:01}. Therefore, by Theorem \ref{FTC}, we conclude that (\ref{FTCeq}) holds for $\mu$ when $\lambda\in\{1-3^{-m}:m\in\mathbb{N}\}$.
\end{example}
\begin{remark}
(1) The authors of  \cite{RS:23} focused on the self-similar measures $\mu$ associated with $(p_i)_{i=1}^3$ and the following IFS:
$g_1(x)=\frac{1}{3}x,\; g_2(x)=\frac{1}{3}x+1,\; g_3(x)=\frac{1}{3}x+3$.
 These measures are exactly re-scaled $\frac{2}{3}$-Cantor measures. Hence, as a particular case of Example \ref{lambdacantor}, (\ref{FTCeq}) holds. (2) By considering $\lambda=3^{-1}$ in (\ref{lambdacantoreq}) with some direct calculations, one can see that, if $E$ is not totally self-similar, the assumption of Remark \ref{z1} (i) may fail and we are not sure whether (\ref{FTCeq}) remains true.
\end{remark}


\end{document}